\pdfoutput=1
\documentclass[12pt, a4paper, openany]{amsart}
\usepackage{amsfonts,amsmath,amssymb,amsthm}
\usepackage{array,longtable,hhline,mathabx}

\usepackage{xcolor}

\usepackage[all]{xypic}
\xyoption{rotate}
\xyoption{pdf}

\usepackage{mathrsfs}
\usepackage{comment}
\usepackage{enumerate}

\theoremstyle{plain}
\newtheorem{theorem}{Theorem}[section]
\newtheorem{cor}[theorem]{Corollary}
\newtheorem{prop}[theorem]{Proposition}
\newtheorem{lemma}[theorem]{Lemma}

\newcounter{thma}
\newtheorem{lproperty}[thma]{Property}

\theoremstyle{definition}
\newtheorem{de}[theorem]{Definition}
\newtheorem{lde}[thma]{Definition} 

\newtheorem*{qdef}{Definition}
\newtheorem{remark}[theorem]{Remark}

\def\A{{\mathscr A}}
\def\B{{\mathscr B}}

\newcommand\GW{{\mathrm{GW}\,}}

\renewcommand{\phi}{\varphi}
\def\m{{\mathfrak m}}

\def\N{{\mathbb N}}
\def\Z{{\mathbb Z}}
\def\H{{\mathcal H}}

\def\F{{\mathfrak F}}

\newcommand{\Ocont}[1]{\mathcal O^{\text{cont}}_{#1}}
\newcommand{\Mcont}[1]{\mathcal M^{\text{cont}}_{#1}}
\newcommand{\Oh}[1]{\mathcal O^{\text{h}}_{#1}}
\newcommand{\Mh}[1]{\mathcal M^{\text{h}}_{#1}}

\newcommand\ch{{\mathrm{char}\ }}
\newcommand\KG{{\mathrm{KG}\,}}
\newcommand\Tot{{\mathrm{Tot\,}\,}}

\newcommand\Spe{{\mathrm{Spec}\,}}
\renewcommand\Im{{\mathrm{Im}\,}}

\title{Grothendieck--Witt Groups of Henselian Valuation Rings}
\date{\today}
\author{Serge Yagunov}

\usepackage{amsmath,pdftexcmds}

\DeclareFontFamily{U}{cbgreek}{}
\DeclareFontShape{U}{cbgreek}{m}{n}{
	<-6>    grmn0500
	<6-7>   grmn0600
	<7-8>   grmn0700
	<8-9>   grmn0800
	<9-10>  grmn0900
	<10-12> grmn1000
	<12-17> grmn1200
	<17->   grmn1728
}{}
\DeclareFontShape{U}{cbgreek}{bx}{n}{
	<-6>    grxn0500
	<6-7>   grxn0600
	<7-8>   grxn0700
	<8-9>   grxn0800
	<9-10>  grxn0900
	<10-12> grxn1000
	<12-17> grxn1200
	<17->   grxn1728
}{}

\DeclareRobustCommand{\Qoppa}{%
	\text{\usefont{U}{cbgreek}{\normalorbold}{n}\symbol{21}}%
}
\makeatletter
\newcommand{\normalorbold}{%
	\ifnum\pdf@strcmp{\math@version}{bold}=\z@ bx\else m\fi
}
\makeatother
\def\Koppa{\Qoppa}

\begin{document}

\begin{abstract}
We show that functors like algebraic $K$-theory (such as unitary or symplectic $K$-functors), as well as the higher Gro\-then\-dieck--Witt groups, possess the local constancy condition for Henselian valuation rings.
Namely, taken with finite coefficients, these functors send canonical residue maps  into isomorphisms. This statement holds in cases of both equal and mixed characteristics.
The proof is based on a slight modification of Suslin's methods. In particular, we use his notion of  universal homotopy.
\end{abstract}

\maketitle

The computation of $K$-groups of various fields has been a classic question in algebraic
$K$-theory. Without delving into the history, let us only mention the work~\cite{quillen1972cohomology} of Daniel Quillen, where he completely computes the $K$-theory of finite fields. 

From an algebro-geometric point of view, $K$-groups of fields carry on some  information about geometric points of an algebraic variety.
The Henselization of a local ring of a point is an infinitesimally small \`etale neighborhood, which is, in most of the cases, non-contractible.  
Thus, besides the $K$-functors of fields, it is also important to study $K$-groups of local rings.
We are trying to answer a natural question whether $K$-functors can
``distinguish" closed points and their infinitesimal neighborhoods. 
For algebraic $K$-theory, the answer to this question was obtained by Andrei Suslin~\cite{suslin1984k}, who has computed, in particular,  the algebraic  $K$-groups with finite coefficients of a local Henselian ring.
As it turned out, these groups coincide with $K$-groups of corresponding residue fields.
So that, taken with finite coefficients, $K$-functors are locally constant. In other words, they are ``non-sensitive'' to the passage from a point to its infinitesimal neighborhood.
Let us write down this property in a more rigorous mathematical form. We reproduce here Definition~\ref{LCC_Functor} below.

\begin{qdef}\label{LCC_Functor}
	We say that the functor $\F$  satisfies the  {\it local constancy condition (LCC)} for a local pair $(R,I)$ if the natural map $R\to R/I$ induces an isomorphism $\F(R)\overset{\cong}{\to} \F(R/I)$.
\end{qdef}

In this paper, we are going to check the LCC for several functors (cohomology theories).
 Among them are  unitary, orthogonal, and symplectic $K$-functors. In these cases, the corresponding computation of $K$-groups for finite
fields was performed by Eric Friedlander~\cite{friedlander1976computations}. 
For these functors, we deduce the LCC from the general rigidity theorem for schemes of finite type (see~ Theorem~\ref{Rigidity} below). We follow the proof strategy of~\cite{suslin1984k}. In particular, we use Suslin's universal homotopy method.  Let us mention here that applying the universal homotopy construction together with the rigidity theorem~\ref{Rigidity}, we obtain a quite easy proof of the LCC for local pairs $(R,I)$ such that $\ch R=\ch R/I$. In the current paper we mostly concentrate at the more subtle case of mixed characteristics, where $\ch R=0$ and $0<\ch R/I\ne m$.

Other functors for which we prove LCC are higher Grothendieck--Witt groups.
In recent years, there has been an increased interest to these functors, due to the results of Fabien Morel relating them to motivic homotopy groups of spheres. We use a general notion of Grothendieck--Witt groups
of exact categories. For them, we deduce the LCC  from classical Suslin's result for $K$-theory, our result about symplectic $K$-groups, and a form of Max Karoubi's induction principle.
In the current paper, we define hermitian $K$-theory via Quillen's ``+''-construction. For Grothendieck--Witt groups we use a construction by Marco Schlichting~\cite{Schlichting2017}.

A brand new approach of~\cite{Calmes2023,Calmes2020,Calmes2020a} allows one to consider Grothendieck--Witt groups in a much broader context.  In~\cite[Definition 1.2.1]{Calmes2023} the hermitian structure is defined as a quadratic functor $\Koppa$ on an $\infty$-category.
In particular, this highly technical approach gives us a mode to treat the case of characteristic two. 
Studying this special case was beyond of our purposes in this article. However, in Remark~\ref{char2} we will touch on a possible strategy for applying mentioned methods to it. 

Recently, utilizing the techniques of~\cite{Calmes2023,Calmes2020,Calmes2020a}, Markus Land~\cite{Land2023} independently obtained similar results. He used the machinery of $\infty$-categories to reduce the LCC property for hermitian $K$-theory to the classical Gabber--Suslin's rigidity theorem for algebraic $K$-theory.
Oppositely, we used more advanced versions of the rigidity theorem~\cite{ananyevskiy2018rigidity,hornbostel2007rigidity,Yagunov2004}
which allowed us to get by with a slight modification of Suslin's approach. 


Finally, let us give a brief outline of the paper.
In the first section we introduce some definitions and conventions, which will be used later.
The second section is technical. There we study homology of algebraic groups over valuation rings and 
check a form of LCC for congruence subgroups. In that section we also restrict ourselves to the case of
Discrete Valuation rings (DVR), the most interesting for the applications. This led to a more clear and self-explanatory proof. However, a similar statement can be established also without this restriction,
analysing filtration spectral sequences.
In the third section, we recall a universal homotopy method and deduce the LCC from the rigidity theorem for 
such functors as $KSp,KU,KO$ (Theorem~\ref{main}). In the last section, we prove our result for the Grothendieck--Witt groups (Theorem~\ref{GW_main}).

\section{Basic notions}

In this section, we assemble together the conventions we use, the necessary properties and general requirements for the objects of our study.

\vskip 5mm

Let  $A(\ast)$ and $B(\ast)$ be two statements parametrized by positive integers.
We say that $A$ {\it  stably implies} $B$ if the assertion: ``$A(\ast)$ holds for all $\ast<N$'' implies
the existence of such a constant $C(N)$ that $B(\ast)$ holds for all $\ast<C(N)$ and $\lim_{N\to\infty}C(N)=+\infty$. Also, if $A$ stably implies $B$ and $B$ stably implies $A$, we say that $A$ and $B$ are {\it  stably equivalent}. 

\subsection*{Rings} All rings are assumed to be unital commutative domains. A pair  $(R,I)$, where the ring $R$ is local and $I$ is an ideal in $R$, will be called a {\it local pair}. If, in addition, the ring $R$ is Henselian and the ideal $I=\mathfrak m$ is maximal, we talk of a Henselian pair.    
For a local ring $R$ we will often denote its field of fractions by $E$ and the residue field $R/\m$ by $F$.

Our main objects are valuation rings (which are automatically local).
Let $(E,v)$ be a field $E$ endowed with a valuation $v$ taking values in $\Gamma\cup\{\infty\}$, where $\Gamma$ is a totally ordered commutative group (with the neutral element 0). 
In this case we say that $R=(R,v)=\{\alpha\in E|v(\alpha)\ge 0\}$ is a valuation ring and denote by $I_\sigma=\{\alpha\in E|v(\alpha)\ge\sigma>0\}$ the valuation ideal in $R$ corresponding to the value $\sigma\in \Gamma$. 
In section~\ref{Homology} we will additionally assume that $R$ is a Discrete Valuation Ring (DVR).    
In this case,  the valuation group $\Gamma\cong\Z$ and the corresponding valuation ideals $I_k=I^k$ are 
powers of the maximal ideal $I_1=I=\mathfrak m$.

\subsection*{Algebraic Groups}

We call $G$ an algebraic group, if $G$ is an affine algebraic ind-group $G=\underset i\lim\, G_i$, having
the following properties.

\begin{itemize}
	\item $G_i$ are affine algebraic groups, which are identified  with algebraic subgroups of $GL_i$. Hence, for a ring $A$ the group $G(A)$ may be identified with the subgroup of $GL(A)$ and $G_i=G\cap GL_i$.   
	\item Groups $GL(A)$ allow the application of Quillen's ``+''-con\-stru\-ction~\cite{kervaire1969smooth,Quillen1994}. In particular, they have perfect commutator subgroups. Applying the ``+''-construction, we obtain sequences of homotopy groups $\pi_\ast(BG^+(A))$ of the spaces $BG^+(A)$, which are called $K$-groups corresponding to $G$. We denote them by $KG_\ast(A)$.
	\item Homology of $G(A)$ (as a discrete group) satisfies the stabilization property below. 
\end{itemize}

\begin{lproperty}\label{Homology_Stability}
	Let $G$ be an algebraic group,  $A$ be a local ring, and $G_i(A):=GL_i(A)\cap G(A)$. Then, for any $n\in\N$ there exists a constant
	$C(n)\in\N$ such that for $0\le i\le n$ and $j>C(n)$ the natural embedding $G_j\subset G$ induces an isomorphism $H_i(G_j(A))\cong H_i(G(A))$.
\end{lproperty}

\subsection*{Theories} We consider functors (cohomology theories) defined at the category of rings (unital commutative domains, affine schemes) and taking values in the category of graded abelian groups.  
Mostly, we will make a deal with  functors $KG_\ast(-):=\pi_\ast(BG^+(-),\Z/m)$ for an appropriate affine algebraic group $G$. 
We consider cohomology theories on the category of algebraic varieties with finite coefficients $\Z/m$. 
Here we list a number of formal properties that our functors must satisfy.  
Let $\mathfrak F$ be such a functor.
Firstly, it is supposed to be rigid. 

\begin{lproperty}[Rigidity Theorem]\label{Rigidity}
	Let $A$ be a commutative unital ring and $\mathfrak F$ be a theory with finite coefficients $\Z/m$ such that $m$ is invertible in $A$.
	Let also $X$ be a smooth scheme of locally finite type over $\Spe A$ and $s\colon \Spe A\to X$ be a section.
	Denote by $X_s^h$ the henselisation of $X$ along $s$. Then, the induced map
	$s^\ast\colon \mathfrak F(X_s^h)\to  \mathfrak F(\Spe A)$
	is an isomorphism.
\end{lproperty}

Secondly, our functor should satisfy a natural analogue of Corollary~\ref{congr_hom_zero}.

For our cohomology theory we are going to check the following property.

\begin{lde}\label{LCC_Functor}
	The functor $\F$  satisfies the  {\it local constancy condition (LCC)} for a local pair $(R,I)$ if the natural map $R\to R/I$ induces an isomorphism $\F(R)\overset{\cong}{\to} \F(R/I)$.
\end{lde}

Most of the known proofs of the rigidity theorem in various cases rely on a number of axioms our functor $\mathfrak F$ must satisfy. As a rule, it is convenient to consider functors that are cohomology theories in the sense of, say~\cite{panin2002rigidity}. They should also satisfy the continuity axiom~\cite[(1.5)]{panin2002rigidity}.

Finally, let us mention that we consider only the theories (functors) with coefficients $\Z/m$. In a typical situation, when the theory $\mathfrak F$ is represented by a space or spectrum $\mathscr F$, we have
$\mathfrak F_\ast(-):=\pi_\ast(\mathscr F(-),\Z/m)$. So, here we need to take homotopy groups with finite coefficients.  Basic information concerning them may be found in~\cite{neisendorfer1980primary}. 
Some discussion of the algebro-geometric case can be also found in~\cite{hornbostel2007rigidity}.
 
\subsection*{Examples}
Currently, our principal example of the cohomology theory 
is the symplectic $K$-theory functor $KSp$.
All the constructions are also applicable to the groups $G=GL,U$, or $O$. 
In all these cases, setting $KG_\ast(R):= \pi_\ast(BG^+(R))$, one gets the corresponding $K$-groups.

It is known that these functors satisfy the list of necessary axioms.

Homological stability theorem was proven in~\cite{nesterenko1990homology,van1980homology} for the group $GL$, in~\cite{mirzaii2002homology} for the case of unitary group, and~\cite{mirzaii2001homology} for symplectic groups. 

The rigidity property for the ``classical'' case $GL$ was proven in~\cite{gabberk}. All other cases 
can be found in~\cite{ananyevskiy2018rigidity,hornbostel2007rigidity}.

\section{On the homology of linear groups}\label{Homology}

In this section, we prove the LCC for homology of linear groups over a discrete valuation ring $R$.    

\begin{de}
	For an affine algebraic group $G$ and a local pair $(R,I)$, let $G(R,I)$ be the kernel of the natural group map
	$G(R)\to G(R/I)$. These (abstract)  groups are called {\it congruence subgroups with respect to the ideal $I$.}
\end{de}

For an algebraic group $G$, we consider homology groups $H_\ast(G(-),\Z/m)$ as
functors $\F_\ast(-)$. We shall find a criterion for functors $\F_\ast$ to satisfy the LCC.   
Here and below we assume that $\frac{1}{m}\in R$.
For brevity, from now on,  we denote the homology groups $H(G(R/I^k),\Z/m)$ by $\H(k)$.

\begin{prop}\label{NilIdeal}
	For every integer $k>0$, the natural map $\phi^k\colon R/I^k\to R/I$ induces an isomorphism in homology $\phi^k_\ast\colon \H(k)\overset\simeq\to \H(1)$. 	
\end{prop}

We start with the following lemma.

\begin{lemma}
	Let  $A$ be a local ring and $J\subset A$ be an ideal. Let us also assume that $J^2=0$.
	 Then, the natural map $H_\ast(G(A))\to H_\ast(G(A/J))$ is an isomorphism.  
\end{lemma}	

\begin{proof}
	Consider the short exact sequence
	$$\xymatrix{1\ar[r]&G(A,J)\ar[r]&G(A)\ar[r]&G(A/J)\ar[r]&1.}$$
	Since $J^2=0$, the congruence subgroup $G(A,J)$ is isomorphic to 
	the additive subgroup of the cotangent vector space $T_e(G(A))$ of the variety $G(A)$ at the unit. This subgroup is an $A/J$-module and hence, it is uniquely $m$-divisible. Therefore, homology groups $\tilde H_\ast(G(A,J))$ with finite coefficients $\Z/m$ vanish and the canonical map is an isomorphism by the Hochschild--Serre spectral sequence argument.  
\end{proof}

\begin{proof}[Proof of the Proposition]
	Applying the previous lemma to the case $A=R/I^{2k}$ and $J=I^k/I^{2k}$,  we conclude that for any integer $k>0$ one has: $\H(2k)\cong \H(k)$.
	
	Let us now show that for any integers $k>n>0$, one has: $\H(k)\cong \H(n)$. Choose an integer $j>0$ such that $2^jk>2^jn>k>n$.
	Then we have the following commutative diagram:
	
	$$\xymatrix{\H(2^jk)\ar@/^2pc/^\alpha_\cong[rr]\ar[r]&\H(2^jn)\ar@/_2pc/_\beta^\cong[rr]\ar^\phi[r]&\H(k)\ar[r]&\H(n).}
	$$ 	
	
	Since $\alpha$ is an isomorphism, we conclude that $\phi$ is an epimorphism. Oppositely, because $\beta$ is an isomorphism, we see that $\phi$ is a monomorphism and  the proposition follows. 
\end{proof}	

\begin{cor}
	If the ideal $I$ is nilpotent, there is an isomorphism 
	$H_\ast(G(R))\overset\cong\to H_\ast(G(R/I))$.
\end{cor}

Now, consider the following. 

\begin{lemma}\label{GroupIsom}
	Denote by $G$ an abstract group, by $V_1,V_2$ its normal subgroups, and by $Q_1,Q_2$ the
	corresponding factor groups.  Assume we are given a morphism  of short exact sequences
	
	$$\xymatrix{1\ar[r]&V_1\ar_\phi@{^(->}[d]\ar[r]&G\ar@{=}[d]\ar[r]&Q_1\ar@{->>}_\psi[d]\ar[r]&1\\
		1\ar[r]&V_2\ar[r]&G\ar[r]&Q_2\ar[r]&1
	}
	$$
	with a monomorphism $\phi$ and an epimorphism $\psi$.
    Let $M$ be an abelian group taken as a trivial $G$-module. Then, the statement i) below stably implies 
    the statement ii).
    \begin{enumerate}[i)]
    	\item  The induced map $$\phi_\ast\colon \tilde H_\ast(V_1,M)\to \tilde H_\ast(V_2,M)$$ is zero and  $$\psi_\ast\colon H_\ast(Q_1,M)\to H_\ast(Q_2,M)$$ is an isomorphism.
    	\item The canonical maps $G\to Q_i$ $(i=1,2)$ induce isomorphisms 
    	$$H_\ast(G,M)\cong H_\ast(Q_i,M).$$
    \end{enumerate}
\end{lemma}

\begin{proof}
	Let $F_\bullet$ be a projective $\Z[G]$-resolution of $M$. Take (for $i=1,2$) the complexes $C^i_\bullet=(F_\bullet)_{V_i}$ of $V_i$-coinvariants.
	The map $\phi$ induces the morphism $\phi_\ast\colon C^1_\bullet\to C^2_\bullet$.
	Denote by $\tilde C^i_\bullet$ a truncation of $C^i_\bullet$ determined as
	
	$$\tilde C^i_n=\begin{cases}
		 C^i_n &\text{ for }n>0;\\
		 \Im d_1&\text{ for }n=0.
	\end{cases}
	$$ 
	One has the following exact sequence of complexes
	$$\xymatrix{0\ar[r]&\tilde{C^i_\bullet}\ar[r]&C^i_\bullet\ar[r]&H_0(C^i_\bullet)\ar[r]&0,}$$
	where the latter term (the group $H_0(C^i_\bullet)=H_0(V_i,M)=M$) is considered as a complex concentrated in the dimension $0$.
	
	Consider homology of groups $Q_i$ with coefficients in these complexes.
	With the help of homomorphisms $\phi_\ast,\psi_\ast$, we obtain a morphism of long exact sequences.
	Consider a fragment of them
	$$\xymatrix{H_\ast(Q_1,\tilde{C^1_\bullet})\ar[r]\ar^{\alpha}[d]&H_\ast(Q_1,C^1_\bullet)\ar^{\zeta}[r]\ar^{\beta}[d]&H_\ast(Q_1,M)\ar[r]\ar^{\gamma}[d]&H_{\ast-1}(Q_1,\tilde{C^1_\bullet})\ar^{\alpha_{-1}}[d]\\
		H_\ast(Q_2,\tilde{C^2_\bullet})\ar[r]&H_\ast(Q_2,C^2_\bullet)\ar[r]&H_\ast(Q_2,M)\ar[r]&H_{\ast-1}(Q_2,\tilde{C^2_\bullet}).
	}
	$$
	Assume, we are in the stable range ($\ast<N$).
	The map $\gamma$ in the diagram is an isomorphism by the lemma premises.
	We also see that $\beta$ is an isomorphism because the groups $H_\ast(Q_i,C^i_\bullet)$ are isomorphic to $H_\ast(G,M)$ (Hochschild--Serre spectral sequence) and
	the map $\beta$ is induced by the identity map on $G$.
	
	Finally, consider the  map 
	$\alpha\colon H_\ast(Q_1,\tilde{C^1_\bullet})\to 
	H_\ast(Q_2,\tilde{C^2_\bullet})$.
	It decomposes as 	
	$$H_\ast(Q_1,\tilde{C^1_\bullet})\overset{\phi_\#}{\to} H_\ast(Q_1,\tilde{C^2_\bullet})\overset{\psi_\#}{\to}	H_\ast(Q_2,\tilde{C^2_\bullet}).$$
	It is not difficult to verify that the first map is stably zero, provided that
	$\phi_\ast\colon \tilde H_\ast(V_1,M)\to \tilde H_\ast(V_2,M)$ equals to zero in the range.
	Then, it results that $\alpha=0$.
	
	We can now complete the proof by the simple diagram chasing, getting that $\zeta$ is an isomorphism.	  
\end{proof}	

Let us get back to the main narrative thread. We have arrived at the following theorem.

\begin{theorem}\label{HomInv}
	Consider a local pair  $(R,I)$  and  an affine algebraic group $G$ over $R$.
	Let for some integers $k>n>0$ the natural embedding $G(R,I^k)\hookrightarrow G(R,I^n)$ induce a
	trivial map on the homology groups $H_j(G(R,I^k))\to H_j(G(R,I^n))$ for $j\le N$. Then, the  induced map  $H_\ast(G(R))\to H_\ast(G(R/I))$ is an isomorphism for $\ast<N$.
\end{theorem}
\begin{proof}
	By Proposition~\ref{NilIdeal}, one has $ H_\ast(G(R/I^k))\cong H_\ast(G(R/I))\cong H_\ast(G(R/I^n))$ and, by the theorem hypothesis, we are in conditions of Lemma~\ref{GroupIsom}.
	The theorem follows from the lemma.
\end{proof}

\section{K-groups of Henselian rings}

We start with the statements connecting $K$-functors and homology of the corresponding linear groups.

\begin{prop}\label{Homotopy_Homology}
	Suppose that $f\colon X\to Y$ be a map of two H-spaces (not necessary preserving the H-structures). Also assume that both the kernel and cokernel of the induced map $\pi_1(X)\to \pi_1(Y)$ are uniquely m-divisible. Then, the following conditions are stably equivalent:
	\begin{itemize}
		\item $f_\ast\colon\pi_\ast(X,\Z/m)\overset{\cong}{\to} \pi_\ast(Y,\Z/m)$;  
		\item $f_\ast\colon H_\ast(X,\Z/m)\overset{\cong}{\to} H_\ast(Y,\Z/m)$. 
	\end{itemize} 
\end{prop}
\begin{proof}
	This is a direct corollary of the Whitehead theorem in the category of topological spaces localized with respect to $\Z/m$. 
\end{proof}	

\begin{cor}\label{congr_hom_zero}

	For a local pair $(R,I)$ and an algebraic group $G$, the following statements are stably equivalent:
	\begin{enumerate}[i)]
		\item $KG_\ast(R)\to KG_\ast(R/I)$ is an isomorphism;
		\item $H_\ast(G(R))\to H_\ast(G(R/I))$ is an isomorphism;
		\item $\tilde H_\ast(G(R,I))=0$.
	\end{enumerate}
As usual, all homotopy and homology groups are taken here with coefficients $\Z/m$.
\end{cor}

\begin{proof}
The $K$-groups are homotopy groups of the spaces $BG^+(R)$. Taking into account that $H_\ast(BG(R))=H_\ast(BG^+(R))$, the stable equivalence $i)\Leftrightarrow ii)$ follows by Proposition~\ref{Homotopy_Homology}.

Hochschild--Serre spectral sequence arguments show that $iii)\Rightarrow ii)$.

Assume now that $H_\ast(G(R))\cong H_\ast(G(R/I))$.
This stably implies $\pi_\ast(BG(R))\cong \pi_\ast(BG(R/I))$.
Since $$BG(R,I)\to BG(R)\to BG(R/I)$$ is a fibration, we get: $\pi_\ast(BG(R,I))$ stably equals to zero.
Hence, again by Proposition~\ref{Homotopy_Homology}, $ii)\Rightarrow iii)$.
\end{proof}

Now we derive our main Theorem~\ref{main} from a statement regarding homology of congruence subgroups of $G$.
Assume that the following proposition holds.

\begin{prop}\label{zero_map}
	Let us fix an arbitrary $N\in \N$. 
	Then, for every $n\in\N$, $0<j<N$, and $\sigma\in \Gamma$ there exists  $\tau\in \Gamma$, $\tau\ge\sigma$ such that the natural embedding
	$G_n(R,I_\tau)\hookrightarrow G(R,I_\sigma)$ induces the zero map in homology
	$$
	\tilde H_j(G_n(R,I_\tau))\overset{0}{\to}\tilde H_j(G(R,I_\sigma)).
	$$
\end{prop}

Then, we can easily pass to the dimension stable case.

\begin{cor}\label{zero_map_stable}
	In the notation above, the natural embedding
	$G(R,I_\tau)\hookrightarrow G(R,I_\sigma)$ induces the zero map in the reduced homology for degrees $<N$.	
\end{cor}
\begin{proof}
	The corollary follows by Theorem~\ref{Homology_Stability}.
\end{proof}

\begin{theorem}\label{main} 
	$\KG_n(R,\Z/m)\cong \KG_n(F,\Z/m)$ for every $n\ge 0$. Here $F=R/\mathfrak m$.
\end{theorem}	

\begin{proof} 
	By Corollary~\ref{zero_map_stable}, we can choose $N>n$ so that for every $\sigma\in\Gamma$ one can find $\tau\ge\sigma$ in such a way that
	$G(R,I_\tau)\hookrightarrow G(R,I_\sigma)$ induces the zero map in homology for degrees $j<n$. Then, by Theorem~\ref{HomInv} we have an isomorphism: $H_j(G(R))\overset{\cong}{\to} H_j(G(R/I))$ and then application of Corollary~\ref{congr_hom_zero} completes the proof of the theorem. 
\end{proof}

So, it remains to prove Proposition~\ref{zero_map}. For this end, we will need some auxiliary concepts.
Let us recall the universal chain homotopy method.  The idea of the construction is, probably, due to Suslin~\cite{suslin1984k}, but we will rather follow a modern exposition by~\cite[page 525 after VI.3.4]{weibel2013k}.

Let $A$ be a ring and $G_n=G\cap GL_n(A)$ be a linear group over it. 
Abusing the notation, we will omit $n$ in the formulas, since it is fixed   
up to Proposition~\ref{homotopy_zero}. In particular, we will write $G$ for $G_n$.

For each integer $i>0$, let us construct the smooth affine scheme $X_i$ over $\Spe A$, setting 
$$
X_{i}:=\underset{i\text{ times}}{\underbrace{G\times\cdots\times G}}.
$$

The scheme $X_i$ is endowed with canonical projections to its $j$-th components $pr_j\colon X_{i}\to X_1^{(j)}$, $(1\le j\le i)$ and the unit section $\Spe A\overset{e}{\to} X_{i}$. There is also the multiplication morphism $\mu\colon X_2\to X_1$
induced by the multiplication in the group $G$.
 
Denote by  $X^h_{i}$ the henselisation of $X_{i}$ with respect to the unit section. 
Both the schemes are affine, so their are the spectra of rings $X_{i}=\Spe\mathcal O_{i}$ and $X^h_{i}=\Spe\Oh{i}$. We also define $\mathcal M^h_{i}$  as the kernel of evident structure map
 $\mathcal O^h_{i}\to A$.
All recently constructed morphisms commute with unit sections. Therefore, one can write similar 
definitions of the projection and multiplication morphisms for schemes $X_i^h$. This enables us to construct
 the desired differential maps below.

Let us  now take an integer $q$ and consider the ring $\mathcal O_{q}$. It contains the elements $\alpha_{i,j}^k$ (for $1\le k\le q,\ 1\le i,j\le n$) 
that are regular functions on $X_q$ corresponding to $(i,j)$ entry of the $k$-th matrix. Consider the universal matrix $\{\alpha_{i,j}\}$ made of the functions $\alpha^1_{i,j}$, corresponding to the variety $X_1$.  Locally, at the unit it is congruent to the identity matrix
modulo the ideal  $\mathcal M^h_{1}$. Therefore, the matrix gives a distinguished element of the congruence group $G(\Oh{1},\Mh{1})$. Obviously, one can also define universal elements $u_i=[\alpha^1,\alpha^2,\ldots,\alpha^i]\in G(\Oh{i},\Mh{i})$, setting $u_i=\prod_{j=1}^i pr^{\ast}_j(\alpha^j)$.

For any discrete group $\mathscr G$, let us denote by $\tilde C_\ast(\mathscr G)$ its standard reduced ($C_0=0$) bar-complex. One has: $\tilde C_k(\mathscr G)=(\Z/m)[\mathscr G^{\times k}]$ and the differential map $d=d^k$ is given by the formula: $d=\sum_{j=0}^k (-1)^jd_j$, where
$$
	d_j(g_1,\ldots,g_k)=\begin{cases} (g_2,\ldots,g_k)&\text{ for }j=0;\\
	                          (g_1,\ldots,g_jg_{j+1},\ldots,g_k)&\text{ for }1\le j<k;\\
	                          (g_1,\ldots,g_{k-1})&\text{ for }j=k.
	                  \end{cases}      
$$

The homology of the complex is the reduced homology $\tilde H_\ast(\mathscr G,\Z/m)$.

Now we are going to turn the bigraded module $$V_{pq}:=\tilde C_p(G(\Oh{-q},\Mh{-q}))$$ into a forth quadrant 
double chain complex. We only need to introduce differential maps $\partial_\ast$ acting in the ``vertical'' direction. Their construction is dual to one of $d_\ast$. By the multiplication morphism $\mu\colon X_2\to X_1$ one obtains the comultiplication operator $\nu\colon \Oh{1}\to \Oh{2}$. (Here we again use the fact that multiplication preserves the unit section.) Similarly, dualizing the formulae for $d^k_j$, one obtains operators $\partial^k_j\colon \Oh{k}\to \Oh{k+1}$. Finally, one should set $\partial^k=(-1)^k\sum_{j=0}^k (-1)^j\partial_j$ The consistency relations $d^2=\partial d+d\partial=\partial^2=0$ are to be checked by the definition.
 
Let us construct the chains $u_i\in V_{i,-i}$, made of $i$ copies of the universal matrix. 
We set $u_i=(\alpha^1,\alpha^2\ldots,\alpha^i)\in V_{i,-i}$ and also define the chain $\mathfrak u=(0,u_1,u_2,\ldots)$, having full degree $0$ in the total complex $\Tot(V)$. We will call this chain $\mathfrak u$ the universal chain.  

Now we shall construct a family of distinguished elements $c_{i}$ lying in $\tilde C_{i+1}(G(\Oh{i},\Mh{i}))$ that serve as contracting homotopy operators for the universal chain $\mathfrak u$.

\begin{lemma}
	Homology groups  $\tilde H_{\ast}(G(\Oh{i},\Mh{i}))$ vanish for all $i>0$.
\end{lemma}
\begin{proof} 
	The lemma is a direct consequence of Property~\ref{Rigidity} and Corollary~\ref{congr_hom_zero}.
\end{proof}

\begin{lemma} The double complex $V_{p,q}$ is acyclic. \end{lemma}
\begin{proof}
By the previous lemma, the double complex $V$ has exact rows. Hence, the double complex spectral sequence converges to zero that proves acyclisity of the whole double complex.	
({Compare~\cite[Proposition VI.3.5]{weibel2013k}}.)
\end{proof}

\begin{cor} The universal chain $\mathfrak u$ is a boundary in the total complex $\Tot(V)$.
	In other words, there exist a family of chains $c_i\in V_{i+1,-i}$ such that for each $i>0$ 
	the equation $u_i=dc_i+\partial c_{i-1}$ holds.	
\end{cor}

\begin{proof}
By the construction, one has the relation $d u_k=\sum_{j=0}^k (-1)^j\partial_j u_{k-1}$.
So that, the universal chain $\mathfrak u$ is a cycle. Then, the corollary follows, since by the statement of the previous lemma the double complex $V$ is acyclic.	
\end{proof}

\begin{remark} 
	Let us mention here that in the equicharacteristic case  our Proposition~\ref{zero_map} immediately follows from the previous construction. Namely, one can show that the natural morphisms 
	$\tilde C_i(G_n(R,I_\tau))\rightarrow \tilde C_i(G(R,I_\sigma))$ are zero-homotopic, using the  elements $c^\prime_{i+1}$ as chain homotopy operators. This implies $\tilde H_\ast(G(R,I))=0$ and, by Corollary~\ref{congr_hom_zero}, the principal result. In this case the proof corresponds almost verbatim to~\cite[Theorem 2.3]{suslin1984k}. 
\end{remark}

From now on, we shall assume that $A$ is a valuation ring. Therefore, one can consider the schemes $X_{\bullet}$ as topological manifolds in the corresponding valuation topology. Denote the local ring of germs of  functions continuous at the unit section $\Spe A\to X_i$  by $\Ocont{i}$, and its maximal ideal by $\Mcont{i}$, as above.

There are canonical morphisms of schemes $X^h_{i}\to X_{i}$ and 
$\Spe\Ocont{i}=X^\mathrm{cont}_{i}\to X_{i}$. Since the ring $\Ocont{i}$ is Henselian 
(see~\cite[VI.\S 4]{raynaud2006anneaux}), one can
decompose every morphism $X^\mathrm{cont}_{i}\to X_{i}$ as $X^\mathrm{cont}_{i}\to X^h_{i}\to X_{i}$. 
Then, we can define the ``universal homotopy elements''  $c^\prime_i\in\tilde C_{i+1}(G(\Ocont{n,i},\Mcont{n,i}))$
as images of $c_{i}$ under the canonical dual map $\Oh{i}\to \Ocont{i}$. 

Now we are ready to prove the following proposition. We return to the initial notation.

\begin{prop}\label{homotopy_zero}
	Let us fix an arbitrary positive integer $N$. Then, for any $n\in\N$ and any $\sigma\in\Gamma$ one can choose such values $\sigma\le\tau\in\Gamma$ and $m\ge n$  that for any $0\le i<N$ the natural embedding
	$$
	\tilde C_i(G_n(R,I_\tau))\hookrightarrow \tilde C_i(G_m(R,I_\sigma)) 
	$$ is null--homotopic. 
\end{prop}

\begin{proof}
Consider the chains $c^\prime_{n,i}$ for $0\le i\le N$. All these chains are germs of functions continuous at 
the unit-sections. Since the congruence subgroups $G_n(R,I_\bullet)$ make a fundamental system of neighborhoods of $e$ in $G_n(R)$, one can find a value $\tau\in \Gamma$ and an integer $m$ such that  all the functions involved are defined in the group $G_n(R,I_\tau)$ and their images lie inside the group $G_m(R,I_\sigma)$.

Let us now consider an elementary chain $\beta$ in $\tilde C_{i}(G_n(R,I_\tau))$ that is an $i$-tuple of
matrices from the congruence group $G_n(R,I_\tau)$ ($R$-rational point of the scheme $X_{n,i}$).
Denote by $c^\prime_{n,i}(\beta)$ the value of the corresponding universal homotopy function at the point $\beta$.
Every chain $b$ in $\tilde C_{i}(G_n(R,I_\tau))$  is some linear combination $b=\sum_j b_j\beta_j$
of elementary chains. The reader can verify that the element $b=\sum_j b_jc^\prime_{n,i}(\beta_j)$ is the desired contracting homotopy operator.
 \end{proof}

Because the  map appearing in the proposition statement is homotopic to zero, Proposition~\ref{zero_map} obviously follows.

\section{Computation of Grothendieck--Witt groups}

In this section all homotopy groups are taken, as before, with finite coefficients. Some terminology and notation below is adopted from~\cite{Schlichting2017}. 

We always denote by $k$ a unital commutative ring
and require that all appearing rings are $k$-modules.
Everywhere, except for the last paragraph, we also assume 
that $\frac{1}{2}\in k$.
 
\begin{de}\label{Admissible}
We call a category $\mathscr A$ {\it admissible} if it is a dg category with weak equivalences and duality such that all sets  $\mathrm{Hom}_{\mathscr A}(-,-)$ are $k$-modules.  In other words, the category $\mathscr A$ is enriched.
In this section, we can always put $k=\Z[\frac{1}{2}]$. 
\end{de}

In~\cite{Schlichting2017}, Schlichting constructed for any admissible category $\mathscr A$  a (topological) spectrum $GW(\mathscr A)$ and defined   Gro\-then\-dieck--Witt groups  of this category with $\Z/m$-coefficients as homotopy groups 
$$GW^{j}_i(\mathscr A,\Z/m):=\pi_i(GW^{[j]}(\mathscr A),\Z/m).$$  
The resulting groups turn out to be 4-periodic: $GW^{[j]}(\mathscr A)\simeq GW^{[j+4]}(\mathscr A)$. 

In~\cite[4.9 and 5.10]{Schlichting2017}, he also constructs a homotopy cartesian square of Grothendieck--Witt spectra with a contractible spectrum in the lower left corner. The square induces the following functorial distinguished triangle in the homotopy category of spectra  

{\small
$$\xymatrix@1{\triangle_\A:\GW(\A)\ar^(.45){I}[r]&\GW(Fun([1],\A))\ar^(.58){\mathrm{Cone}}[r]&\GW^{[1]}(\A)
	\ar[r]&S^1\land\GW(\A).}
$$}

By~\cite[6.2]{Schlichting2017}, the second term $\GW(Fun([1],\A))$ of $\triangle_\A$ may be identified  with the K-theory spectrum $K(\A)$. After this identification, the morphisms $I$ and $\mathrm{Cone}$
in the triangle become the forgetful and hyperbolic dg form functors $F,H$ (see~\cite[Sect. 4.7]{Schlichting2017})
	$$\xymatrix@1{\GW(\A)\ar^{F}[r]&K(\A)\ar^(.4){H}[r]&\GW^{[1]}(\A).}$$

It is just left to describe the morphism $\GW^{[1]}(\A)\to S^1\land\GW(\A)$ in our triangle. 
Letting $\A$ be the shifted category of bounded complexes of finitely generated free $k$-modules (comp. $\A=\mathscr C^{[-1]}_k$ in~\cite[Ex. 1.2]{Schlichting2017}), one defines
the Bott element  as $\eta:=-\delta(\langle 1\rangle)\in\GW^{-1}_{-1}(k)$. Here  
$\langle 1\rangle\in\GW^{0}_{0}(k)$  is the Grothendieck--Witt class  represented by the symmetric
form $x \otimes y \to xy$ on the ring $k$ and $\delta$ is the boundary map in the fibration exact sequence
corresponding to the cartesian square of {\it loc.cit.}~4.9.
Finally, the desired map is induced by $\cup \eta$.

We can summarize all our observations in the following theorem (Comp.~\cite[Thm. 6.1]{Schlichting2017}).

\begin{theorem}[Algebraic Bott Sequence]\label{ABS} Let $\A$ be an admissible category. Then the sequence of spectra
	$$\xymatrix@1{\triangle_\A\colon\GW^{[n]}(\A)\ar^(.6){F}[r]&K(\A)\ar^(.4){H}[r]&\GW^{[n+1]}(\A)\ar^(.45){\eta\cup}[r]&S^1\land\GW^{[n]}(\A).}
	$$ 
	is an exact functorial triangle in the homotopy category of topological spectra.
\end{theorem}

For a commutative unital ring $R$  denote by $\mathscr P(R)$ the category of finitely generated projective modules over $R$. This category has a standard duality functor $\mathrm{Hom}_R(-,R)$ and  satisfies the conditions of Definition~\ref{Admissible}.	
Morphisms $F$ and $H$ are  natural transformations, by the construction. Let us check the functoriality  of the Bott element.  It suffices to mention that for the ring homomorphism $f\colon R\to S$, one has: 
$f_\ast(\langle 1_R\rangle)=\langle 1_S\rangle$ in the ``classical'' Grothendieck--Witt ring $\GW^0(S)$. 
So that, the previous theorem implies the following corollary.

\begin{cor} The ring homomorphism  $f\colon R\to S$ induces the morphism 
	$f_\#\colon \triangle_{\mathscr P(R)}\to \triangle_{\mathscr P(S)}$ of distinguished triangles of spectra and, therefore, the morphism of two long exact sequences of Grothendieck--Witt groups.
\end{cor}

\begin{theorem}\label{GW_main} 
	Let $(R,I)$ be a Henselian pair such that $\frac{1}{2},\frac{1}{m}\in R$ and
	the ideal $I$ is maximal. Then, for $j\ge 0$ the  functors $\GW_j^\ast(-,\Z/m)$ satisfy the local constancy condition (LCC). 
\end{theorem}	 

In other words, the theorem tells us the following.

\begin{cor}
Let $(R,I)$ satisfy the conditions of the previous theorem. Then, for any integers $j\ge 0$ and $n$ there are isomorphisms $\GW_j^{n}(R,\Z/m)\cong \GW_j^{n}(F,\Z/m)$ induced by the natural map $R\to F=R/I$. 
\end{cor}

The theorem is an easy consequence of the following lemmata. 
The first lemma is a case of Karoubi's induction statement.
\begin{lemma}\label{induction}
	Let $\A, \B$ be two admissible categories and  $\mathfrak F\colon\A\to\B$ be an exact dg form functor.
	Assume that $\mathfrak F$ induces an equivalence of $K$-theory spectra and isomorphisms
	$\GW^{[n]}_{\ge 0}(\A) \cong \GW^{[n]}_{\ge 0}(\B)$ for some $n$.
	Then, $\mathfrak F$ induces the corresponding  isomorphisms for all $n$. 
\end{lemma}
\begin{proof} Let us write down a morphism $\triangle_\A\to\triangle_\B$ of distinguished triangles from Theorem~\ref{ABS}.
	
	$$\xymatrix{\GW^{[n-1]}(\A)\ar^{}[r]\ar[d]&K(\A)\ar^{}[r]\ar^[@!-90]\simeq[d]&\GW^{[n]}(\A)\ar^{}[r]\ar[d]&S^1\land\GW^{[n-1]}(\A)\ar[d]\\
	\GW^{[n-1]}(\B)\ar^{}[r]&K(\B)\ar^{}[r]&\GW^{[n]}(\B)\ar^{}[r]&S^1\land\GW^{[n-1]}(\B).}
$$ 

Passing to the homotopy groups, we obtain a morphism of long exact sequences.
 Obviously, if the vertical arrows induce isomorphisms on the groups $\GW^n_{i+1}$ and $\GW^n_{i}$, then the same is true for the term $\GW^{n-1}_{i}$. 
	By induction and the 4-periodicity of $\GW^{[n]}(-)$, we get the desired statement.
\end{proof}

Our last lemma is borrowed from~\cite[Corollary A.2]{Schlichting2017}.

\begin{lemma}\label{KSpKO}
	For any ring $R$, satisfying our conditions and an integer $j\ge 0$ there are natural isomorphisms $KSp_j(R)\cong \GW_j^{2}(R)$ and $KO_j(R)\cong \GW_j^{0}(R)$.
\end{lemma}

\begin{proof}[Proof of Theorem~\ref{GW_main}.]
	Denote, as usual, our Henselian local ring by $R$ and its residue field by $F$. 
	We can use either $KSp$, or $KO$ theory for our proof. Let us assume, we have chosen $KSp$.
	By Suslin~\cite{suslin1984k}, the map $R\to F$ induces an equivalence of K-theory spectra $K(R)$ and $K(F)$.
	By Theorem~\ref{main}, one has natural isomorphisms $KSp_i(R)\cong KSp_i(F)$ for all $i\ge 0$.
	 Then, by Lemma~\ref{KSpKO}, we are in the conditions of Lemma~\ref{induction} for the case $n=2$. Hence, the theorem follows.
\end{proof}

\begin{remark}\label{char2}
Here we  briefly discuss what happens in the case $\ch F=2$.
There we can no longer use most of the results,
on which we have relied in this chapter. Instead, we can apply the new approach to the Hermitian $K$-theory, recently developed 
in~\cite{Calmes2023,Calmes2020,Calmes2020a}.
We give an outline of the proof strategy.

Let us start from the symplectic $K$-theory $KSp(R)$ which we defined using the ``+''-construction. First of all, by classical methods, we can show the equivalence of this construction to the construction based on Karoubi's~\cite{Karoubi1980} 
groupoids and, therefore, to the ``classical'' Grothendieck--Witt spaces
$\mathcal{GW}^{\lambda}_\mathit{cl}$ (see \cite[p.2]{Hebestreit2021}). Then, one can use the equivalence of \cite[Theorem A]{Hebestreit2021} and pass to spaces $\mathcal{GW}(\mathcal{D}^p(R),\Koppa^{g\lambda})$. 

Now, instead of a morphism $\triangle_\mathscr{A}\to \triangle_\mathscr{B}$
of the exact triangles from Lemma~\ref{induction}, one obtains a morphism
of corresponding Bott--Genauer sequences~\cite{Calmes2020}. 
This enables us to write an induction lemma similar to~\ref{induction}.
Finally, the quasi-periodicity result~\cite[3.5.17]{Calmes2023} proves the LCC property for Poincare structures $\Koppa^{\ge m}$. Passing to the (co-)limit
extends the property to the symmetric and quadratic Poincare structures.
\end{remark}

\section*{Acknowledgement} I am deeply indebted to Baptiste Calm\`es who 
introduced me to the technique of~\cite{Calmes2023,Calmes2020,Calmes2020a} and its applicability to the case of characteristic two. I also want to thank Ivan Panin and Andrei Lavrenov for numerous discussions and valuable comments.

\bibliographystyle{unsrt}

\bibliography{K_theory}

\begin{thebibliography}{10}

\bibitem{quillen1972cohomology}
Daniel Quillen.
\newblock On the cohomology and {$K$}-theory of the general linear groups over
  a finite field.
\newblock {\em Annals of Mathematics}, 96(3):552--586, 1972.

\bibitem{suslin1984k}
Andrei~A Suslin.
\newblock On the {$K$}-theory of local fields.
\newblock {\em Journal of pure and applied algebra}, 34(2-3):301--318, 1984.

\bibitem{friedlander1976computations}
Eric~M Friedlander.
\newblock Computations of {$K$}-theories of finite fields.
\newblock {\em Topology}, 15(1):87--109, 1976.

\bibitem{Schlichting2017}
Marco Schlichting.
\newblock Hermitian {$K$}-theory, derived equivalences and
  {K}aroubi{\textquotesingle}s fundamental theorem.
\newblock {\em Journal of Pure and Applied Algebra}, 221(7):1729--1844, jul
  2017.

\bibitem{Calmes2023}
Baptiste Calm{\`e}s, Emanuele Dotto, Yonatan Harpaz, Fabian Hebestreit, Markus
  Land, Kristian Moi, Denis Nardin, Thomas Nikolaus, and Wolfgang Steimle.
\newblock Hermitian {$K$}-theory for stable $\infty$-categories {I}:
  Foundations.
\newblock {\em Selecta Mathematica}, 29(1):10, 2023.

\bibitem{Calmes2020}
Baptiste Calm{\`e}s, Emanuele Dotto, Yonatan Harpaz, Fabian Hebestreit, Markus
  Land, Kristian Moi, Denis Nardin, Thomas Nikolaus, and Wolfgang Steimle.
\newblock Hermitian {$K$}-theory for stable $\infty$-categories {II}: Cobordism
  categories and additivity.
\newblock {\em arXiv preprint arXiv:2009.07224}, 2020.

\bibitem{Calmes2020a}
Baptiste Calmes, Emanuele Dotto, Yonatan Harpaz, Fabian Hebestreit, Markus
  Land, Kristian Moi, Denis Nardin, Thomas Nikolaus, and Wolfgang Steimle.
\newblock Hermitian {$K$}-theory for stable $\infty$-categories {III}:
  Grothendieck-{W}itt groups of rings.
\newblock {\em arXiv preprint arXiv:2009.07225}, 2(3), 2020.

\bibitem{Land2023}
Markus Land.
\newblock Gabber rigidity in hermitian {$K$}-theory.
\newblock {\em Proceedings of the Royal Society of Edinburgh: Section A
  Mathematics}, pages 1--6, 04 2023.

\bibitem{ananyevskiy2018rigidity}
Alexey Ananyevskiy and Andrei Druzhinin.
\newblock Rigidity for linear framed presheaves and generalized motivic
  cohomology theories.
\newblock {\em Advances in Mathematics}, 333:423--462, 2018.

\bibitem{hornbostel2007rigidity}
Jens Hornbostel and Serge Yagunov.
\newblock Rigidity for henselian local rings and {$\mathbb A^1$}-representable
  theories.
\newblock {\em Mathematische Zeitschrift}, 255(2):437--449, 2007.

\bibitem{Yagunov2004}
Serge Yagunov.
\newblock Rigidity {II}: non-orientable case.
\newblock {\em Documenta Mathematica}, 9:29--40, 2004.

\bibitem{kervaire1969smooth}
Michel~A Kervaire.
\newblock Smooth homology spheres and their fundamental groups.
\newblock {\em Transactions of the American Mathematical Society}, 144:67--72,
  1969.

\bibitem{Quillen1994}
Daniel Quillen.
\newblock On the group completion of a simplicial monoid.
\newblock {\em preprint}, 1994.

\bibitem{panin2002rigidity}
Ivan Panin and Serge Yagunov.
\newblock Rigidity for orientable functors.
\newblock {\em Journal of Pure and Applied Algebra}, 172(1):49--77, 2002.

\bibitem{neisendorfer1980primary}
Joseph Neisendorfer.
\newblock {\em Primary homotopy theory}, volume 232.
\newblock American Mathematical Soc., 1980.

\bibitem{nesterenko1990homology}
Yu~P Nesterenko and Andrei~A Suslin.
\newblock Homology of the full linear group over a local ring, and {M}ilnor's
  {$K$}-theory.
\newblock {\em Mathematics of the USSR-Izvestiya}, 34(1):121--145, 1990.

\bibitem{van1980homology}
Wilberd van~der Kallen.
\newblock Homology stability for linear groups.
\newblock {\em Inventiones mathematicae}, 60(3):269--295, 1980.

\bibitem{mirzaii2002homology}
Behrooz Mirzaii and Wilberd van~der Kallen.
\newblock Homology stability for unitary groups.
\newblock {\em Doc. Math}, 7:143--166, 2002.

\bibitem{mirzaii2001homology}
Behrooz Mirzaii and W~Van~der Kallen.
\newblock Homology stability for symplectic groups.
\newblock {\em arXiv preprint math/0110163}, 2001.

\bibitem{gabberk}
Ofer Gabber.
\newblock K-theory of henselian local rings and henselian pairs.
\newblock {\em Contemp. Math.}, 126:59--70, 1992.

\bibitem{weibel2013k}
Charles~A Weibel.
\newblock {\em The {$K$}-book: An introduction to algebraic {$K$}-theory},
  volume 145.
\newblock American Mathematical Society Providence, 2013.

\bibitem{raynaud2006anneaux}
Michel Raynaud.
\newblock {\em Anneaux locaux hens{\'e}liens}, volume 169.
\newblock Springer, 2006.

\bibitem{Karoubi1980}
Max Karoubi.
\newblock Le th{\'e}oreme fondamental de la {$K$}-th{\'e}orie hermitienne.
\newblock {\em Annals of mathematics}, 112(2):259--282, 1980.

\bibitem{Hebestreit2021}
Fabian Hebestreit and Wolfgang Steimle.
\newblock Stable moduli spaces of hermitian forms.
\newblock {\em arXiv preprint arXiv:2103.13911}, 2021.

\end{thebibliography}

\end{document}